\documentclass[mathpazo]{cicp}
\usepackage{color}
\usepackage{amsmath}
\usepackage{amssymb}
\usepackage{amsthm}
\usepackage{amsfonts}
\usepackage{graphicx}

\usepackage{algorithm}
\usepackage{subfigure}
\usepackage{units}
\usepackage{multirow}
\usepackage{graphicx, epsfig}

\usepackage{bm}
\usepackage[colorlinks=true]{hyperref}
\newtheorem{thm}{Theorem}[section]

\newtheorem{lem}[thm]{Lemma}

\theoremstyle{definition}

\numberwithin{equation}{section}
\def\bsA{\boldsymbol{A}}

\newcommand{\ol}[1]{%
\mbox{$\overline{#1}$}}

\providecommand{\wv}{}
\providecommand{\hbz}{}
\providecommand{\edt}{}
\begin{document}
\title[Iterative solution of Sch\"odinger equation]{A preconditioned iterative solver for the scattering solutions of the Schr\"odinger equation}


\author[Hisham bin Zubair et.~al.]{Hisham bin Zubair\affil{1}\comma\corrauth, Bram Reps\affil{2}, and Wim Vanroose\affil{2}}
\address{\affilnum{1}\ Department of Mathematical Sciences, Faculty of Computer Science, Institute of Business Administration, University Rd., 75270 Karachi, Pakistan. \\
\affilnum{2}\ Department of Mathematics and Computer Science, Universiteit Antwerpen, Middelheimlaan 1, B-2020 Antwerpen, Belgium.}
\emails{\\{\tt h.binzubair@gmail.com} (H.~bin Zubair) \\ {\tt bram.reps@ua.ac.be} (B.~Reps) \\ {\tt wim.vanroose@ua.ac.be} (W.~Vanroose)}


\begin{abstract}
The Schr\"odinger equation defines the dynamics of quantum particles which
 has been an area of unabated interest in physics. We demonstrate how simple 
transformations of the Schr\"odinger equation leads to a coupled linear system,
 whereby each diagonal block is a high frequency Helmholtz problem. Based on this model, 
we derive indefinite Helmholtz model problems \hbz{with strongly varying wavenumbers. We employ the iterative approach for their solution. In particular, we develop} a  preconditioner that has its spectrum restricted to a quadrant \hbz{(of the complex plane) thereby making 
it easily invertible by multigrid methods with standard components. This multigrid preconditioner is used in conjuction with suitable Krylov-subspace methods for solving the indefinite Helmholtz model problems. The aim of this study is to report the feasbility of this preconditioner for the model problems.} We compare this idea with the other prevalent preconditioning ideas, \hbz{and discuss its merits}. Results of numerical experiments are presented, which complement the proposed ideas, and show that this preconditioner may be used in an automatic setting.
\end{abstract}


\maketitle

\section{Introduction}
\label{sec1}
Acoustic, electromagnetic or seismic waves can all be modeled by a
Helmholtz equation with a wave number that has properties specific to
the problem area.  In some acoustic scattering applications, for
example, the wave number is space independent but the boundary of the
domain can be very complicated depending on the shape of the object.
In electromagnetic scattering there are jumps in the material
parameters which lead to a piecewise constant wavenumber.  In a
similar way, the wavenumber in seismic waves will contain information
about the geological layers in the earths crust.  Each of these
problems pose different challenges to the numerical methods.

In this article, we focus on the iterative solution of the Helmholtz
equations with a wave number that is specific to models for breakup
problems in chemical systems. These breakup dynamics are described by
a Schr\"odinger equation that reduces, in the energy range of breakup
problems, to Helmholtz equation with a wavenumber that is continuous
in the space variables and can become very large near the boundary of
the domain. One example is the disintegration into four charged
particles of the H$_2$ molecule when it is hit with a single photon
\cite{Wim05}.

The \hbz{prevalent practice for solving this type of problem} requires massively parallel
computers \cite{taylor2002computational} and \hbz{they use}
a significant portion of the resources of large computer
facilities.  The long term aim of \hbz{this research} is
to replace this practice by efficient iterative methods.

The Helmholtz equation \hbz{has often} outgoing wave boundary conditions.
Fixing homogeneous Dirichlet boundary conditions, on the boundaries of
the truncated numerical domain, leads to artificial reflections in the
domain of interest. These reflections are numerical errors and must be
diminished by damping the outgoing waves at the domain boundaries. To
bring this about in our numerical solution method, we employ {\it
  exterior complex scaled} \cite{S79} absorbing boundary layers
(henceforth ECS-ABL). There is a long history of this type of
absorbing boundary condition for chemical reactions
\cite{moiseyev1998}. This treatment is equivalent to the use of
perfectly matched layers (PML) \cite{B94,CW94} and leads to a
non-Hermitian discrete problem \cite{reps2009}.  For a review on
transparant and absorbing boundary conditions for the Schr\"odinger
equaton we refer to \cite{antoine2008review}.
  
For \hbz{Krylov-subspace methods, the main challenge is to} find a good preconditioner.
Over the years there have been different approaches to preconditioning the \hbz{indefinite} Helmholtz equation. One line of research is based on a shifted Laplacian preconditioner that started with the work \cite{Bayliss83, Bayliss85} ( Bayliss, Goldstein and Turkel).  
They used the Laplacian and the positively shifted Laplacians as preconditioner. 

This was later successfully generalized into a robust method, known as
the {\it complex shifted laplacian} ($CSL$) preconditioner, by Erlangga,
Vuik and Oosterlee using complex valued shifts in
\cite{Yogi04,Yogi06}. Introducing a complex shift pushes the spectrum
of the Helmholtz operator into a region \hbz{that is favorable} for multigrid methods \cite{Achi77,Stu82,Trot01} to
approximately invert the preconditioning problem. It is well-known that multigrid efficiency can readily be exploited \hbz{only} for problems having (positive or negative) definite \hbz{spectra}. In the indefinite case \cite{Trot01}, both vital components of multigrid, i.e., smoothing, and coarse grid correction suffer severe degradation, and consequently this results in
divergence of the method \cite{EEL01}.

An alternative preconditioner that, in addition to shift, also scales
the Laplacian was derived from frequency shift time integration by
Meerbergen and Coyette \cite{meerbergen2009}.  By appropriately
choosing the shift and the scale it is possible to restrict the
spectrum of preconditioning matrix into one quadrant of the complex
plane.  We call this type of preconditioner a {\it quadrant definite}
($QD$) preconditioner.

In \cite{reps2009}, we proposed the {\it complex-scaled grid}\, ($CSG$) preconditioner, 
and demonstrated its utility in connection with indefinite Helmholtz 
problems constructed with ECS-ABL. Both the $CSL$ and the $CSG$ preconditioners have similar performance 
and are based on similar ideas. The $CSL$ \textit{translates} the spectrum, while the $CSG$ \textit{rotates}
 it, \hbz{thereby placing it in a region which is multigrid favorable}. Both of these preconditioners depend on the 
translation magnitude or the rotation angle which has to be tuned for specific problems. 

This paper \wv{studies} a preconditioner based on a scaled translation
of the spectrum \wv{that} restricts it to one quadrant of the complex
plane. \wv{We evaluate it on a set of model problems representative for
breakup problems that are derived in the paper.}  \hbz{While its efficiency is found to be between that of the Laplacian
preconditioner and the $CSL$/$CSG$ preconditioners, the main merit is its ease of invertibility by multigrid methods that use well-known standard components. This is a clear advantage of using the $QD$ method, as for the $CSL$/$CSG$ preconditioners, multigrid has to be tuned for different wavenumbers. Moreover, a shift for the $CSL$ preconditioner (or equivalently, a rotation angle for the $CSG$ preconditioner) is apparantly only available through a hit-and-trial rule. In comparison, the $QD$ preconditioner may be used in an automatic setting.}

\wv{Both the discretisation and the absorbing boundary conditions used
  in this paper are of low order of accuracy. Both can be replaced by
  higher order methods, however, the focus of the paper is on the
  working of the iterative methods and this can be studied with the
  low order methods since the higher order discretisation and
  boundary conditions have similar spectral properties.}

In Section \ref{sec:derive}, we give the transformation of the
Schr\"odinger equation to a coupled Helmholtz problem, and derive the
model problems for this study. The details of ECS-ABL and \hbz{the discretization} are given in Section \ref{sec:discretize}. Also reviewed
here, are the spectral properties of the discrete operator. Next, in
Section \ref{sec:precmg}, we describe the $QD$ preconditioner in detail,
and give the multigrid algorithm which we use for approximate
inversion of the preconditioners. This is followed by numerical
experiments, which are given in Section \ref{sec:numex}. Some
conclusions mark the end of the paper in Section \ref{sec:conclude}.
 
\section{From the Schr\"odinger equation to a coupled Helmholtz problem}
\subsection{The model problem}
\label{sec:derive}
In this section we derive the model Helmholtz problem that we use in
this paper to benchmark iterative solvers.  The model problem is \wv{
\begin{equation}\label{eq:modelproblem}
\left\{\begin{split}
 \left(-\Delta_{l_1,l_2} - k^2(x,y) \right) u(x,y) = f(x,y) \quad \text{on}\quad [0,a]^2 \subset \mathbb{R}^2,\\
 u(x,0) = 0 \quad  \forall x \in [0,a]  \quad \text{and}\quad  u(0,y) = 0 \quad \forall y \in [0,a],\\
 \text{ABC on} \quad u(x,a)\quad   \forall x \in [0,a]  \quad \text{and} \quad  u(a,y)\quad   \forall y \in [0,a],
 \end{split}
\right.
\end{equation}
where  ABC denotes outgoing wave boundary conditions, see Section \ref{sec:discretize}, and 
}
\begin{equation}\label{eq:definition_delta}
 \Delta_{l_1,l_2} = \partial_{xx} + \partial_{yy} -\frac{l_1(l_1+1)}{ x^2} - \frac{l_2(l_2+1)}{y^2}
\end{equation}
\wv{ 
denotes the radial part of the Laplacian in spherical coordinates with $l_1$, $l_2 \in \mathbb{N}$.
The wavenumber $k^2(x,y) = 2m(E-V(x,y))$ depends on a potential
$V(x,y)$ that} \hbz{varies continuously in} \wv{$x$ and $y$}
\hbz{in the domain }\wv{ $[0,a]^2$,  $E>0$ is the energy and $m>0$ is the mass of the system.  The right hand side $f(x,y)$ is assumed
to be zero outside $[0,b]^2$ with $b<a$ so that the Helmholtz problem becomes a
 homogeneous problem  in a strip near the boundaries with the ABC.}
%
%
\subsection{The Schr\"odinger equation}
To derive this model we start from the driven Schr\"odinger equation
\begin{equation}\label{eq:schrodinger}
 (\mathcal{H}-E)\psi(\mathbold{r}_1, \mathbold{r}_2) =  \phi(\mathbold{r}_1,\mathbold{r}_2),
\end{equation}
\wv{with $\mathbold{r}_1$, $\mathbold{r}_2 \in \mathbb{R}^3$} and where $\mathcal{H}$ denotes the Hamiltonian and is given by
\begin{equation}
  \mathcal{H} =  - \frac{1}{2m} \Delta_{\mathbold{r}_1}  - \frac{1}{2m} \Delta_{\mathbold{r}_2}  + V_1(|\mathbold{r}_1|) + V_2(|\mathbold{r}_2|) + V_{12}(\mathbold{r}_1,\mathbold{r}_2)
\end{equation}
with $V_1$ and $V_2$ local potentials that only depend on magnitude of
$\mathbold{r}_i$. The potential $V_{12}$ depends, \wv{usually,} on the relative
distance between $\mathbold{r}_1$ and $\mathbold{r}_2$. \wv{The mass $m>0$ scales the Laplacians.} The right hand
side of \eqref{eq:schrodinger}, \wv{$\phi(\mathbf{r}_1,\mathbf{r}_2)$,
 is assumed zero if $|\mathbf{r}_1| >b$ or $|\mathbf{r}_2| > b$,} and can
model an incoming electron that impacts in the system \cite{rescigno1999}, \hbz{or alternately, represents the dipole operator working on a groundstate if the model is used to compute photo-ionization \cite{Wim05}}. 

For these breakup problems the solution
$\psi(\mathbold{r}_1,\mathbold{r}_2)$ is  an outgoing wave in any
direction similar to the Sommerfeld radiation condition.  This leads to a six dimensional problem on an unbounded
domain.  The problem can also be interpreted as a 6D Helmholtz problem
\begin{equation}
  \left(-\Delta_{6D} - k^2(\mathbold{r}_1,\mathbold{r}_1)\right)\psi(\mathbf{r}_1,\mathbf{r}_2) = f(\mathbold{r}_1,\mathbold{r}_2),
\end{equation}
where $k(\mathbold{r}_1,\mathbold{r}_2) = \sqrt{2m(E-V)}$ with $V$
denotes the sum of all potentials.  This becomes a Helmholtz problem
with a constant wave number, $k=\sqrt{2mE}$, in the regions of space
where the potentials go to zero. This  6D problem is hard to solve with the current generation of numerical methods.

%
%
\subsection{Expansion of the solution in partial waves.}
In this section we discuss the reduction of the 6D problem to a
coupled set of 2D problems.  At large distances \hbz{the} solution behaves as
a spherical \hbz{wave} emerging from the center of mass of the system.  It
is therefore common practice \cite{baertschy2001,vanroose2006double}
to rewrite the equation \eqref{eq:schrodinger} in spherical
coordinates.  \wv{The Laplacian operator then splits into a radial
operator and the angular operator differential  \cite{arfken}.}  The coordinates are
written as $\mathbold{r}_1 = (\rho_1, \Omega_1)$ and $\mathbold{r}_2 =
(\rho_2, \Omega_2)$, where $\Omega$ denotes $(\theta, \varphi)$.
\wv{The solution is then written} as a series
\begin{equation}\label{eq:proposal}
 \psi(\mathbold{r}_1,\mathbold{r}_2) = \wv{\sum_{l_1=0}^{\infty}\, \sum_{m_1=-l_1}^{l_1}\,\sum_{l_2=0}^\infty\, \sum_{m_2=-l_2}^{l_2} }  \psi_{l_1m_1,l_2m_2}(\rho_1,\rho_2) Y_{l_1m_1}(\Omega_1) Y_{l_1m_1}(\Omega_2), 
\end{equation}
where $Y_{lm}(\Omega)$ are the spherical harmonics, \wv{the
  eigenfunctions of the angular differential operator of the Laplacian
  in spherical coordinates \cite{arfken}.  In physics this decomposition is referred
  to as the partial wave expansion and the functions
  $\psi_{l_1\,m_1,l_2\,m_2}(\rho_1,\rho_2)$ are called partial waves.}

When this proposal, \eqref{eq:proposal}, is \hbz{substituted} in \eqref{eq:schrodinger} we find
an equation for $\psi_{l_1\,m_1,l_2\,m_2}(\rho_1,\rho_2)$ for all
$l_1 \ge 0,l_2 \ge 0 ,|m_1| \le l_1 ,|m_2| \le l_2$ that is
coupled to all other partial waves.
\begin{equation}\label{eq:coupled}
\begin{aligned}
 &\left[-\frac{1}{2m}  \Delta_{l_1,l_2} -E\right]  \psi_{l_1 m_1,l_2 m_2}(\rho_1,\rho_2) \\
&+ \sum_{l^{\prime}_1=0}^{\infty} \sum_{m_1=-l_1}^{l_1} \sum_{l_2^\prime = 0}^\infty \sum_{m^\prime_2=-l_2^\prime}^{l_2^\prime} V_{l_1\,m_1\,l_2\,m_2 ; l_1^\prime\,m_1^\prime,l_2^\prime\, m_2^\prime}(\rho_1,\rho_2)\,\, \psi_{l_1^\prime m_1^\prime,l_2^\prime m_2^\prime}(\rho_1,\rho_2)  = \varphi_{l_1m_1,l_2m_2},
\end{aligned}
\end{equation}
where the coupling potentials are calculated as
\begin{equation}\label{eq:coupled_potential}
\begin{aligned}
 V_{l_1m_1l_2 m_2;l_1^\prime m_1^\prime l_2^\prime m_2^\prime}(\rho_1,\rho_2 )= &\int d\Omega_1 d\Omega_2  Y^*_{l_1m_1}(\Omega_1)Y^*_{l_2m_2}(\Omega_2)\\
&\times \left[V_1(|\mathbold{r}_1|) + V_2(|\mathbold{r}_2|) + V_{12}(\mathbold{r}_1,\mathbold{r}_2)\right]Y_{l^\prime_1 m^\prime_1}(\Omega_1)Y_{l_2m_2}(\Omega_2) 
\end{aligned}
\end{equation}
and $\varphi_{l_1\,m_1,l_1\,m_2}$ is partial wave of the right hand side.

When the potentials $V_1$, $V_2$ and $V_{12}$ are spherical\hbz{ly} symmetric
the system decouples. When it is c\hbz{y}lindrical\hbz{ly} symmetric the different
$m_1$ and $m_2$ are decoupled. \wv{But in general the system is fully coupled.} Furthermore, it is common practice to
truncate the infinite series in $l$ at a finite $l_\text{max}$ so that it
becomes a finite system of coupled partial differential equations.

\wv{ The boundary conditions \hbz{for Equation} \eqref{eq:schrodinger} translate in spherical coordinates into
  homogeneous Dirichlet $\psi(\rho_1,0)=0$ for all $\rho_1$ and
  $\psi(0,\rho_2)=0$ for all $\rho_2$.  \wv{This is typical for radial problem since $rh0_1=0$ and $\rho_2=0$ is now the origin of the coordinate system \cite{arfken}.} The outgoing boundary conditions
  translate then into outgoing boundary conditions $\rho_1 \rightarrow
  \infty$ or $\rho_2 \rightarrow \infty$. We will \hbz{elaborate on this
  topic in Section \ref{sec:discretize}}.}

\wv{The partial wave expansion can also be written down for a single
  particle Hamiltonian. It then involves an expansion over \hbz{a} single
  angular function $Y_{lm}$ \hbz{and subsequently} leads to a coupled system of
  ordinary differential equations. On the other hand, the Hamiltonians
   currently studied in the physics and chemistry communit\hbz{ies involve}
  three or more particles.  For three particles the driven
  Schr\"odinger equation is a 9-dimensional equation that, after
  expansion in partial waves, becomes a set of coupled 3D PDEs.}

\subsection{Blocked structure and Iterative methods}
The system \eqref{eq:coupled} has a very particular structure.  Since
the differential operators are block diagonal in the spherical
expansion, they only appear on the diagonal blocks of the
equation. The blocks are only coupled by the potentials \wv{defined in
  equation \eqref{eq:coupled_potential}}. \wv{The Hamiltonian
  $\mathcal{H}$ can be written in blocked matrix notation as
\begin{equation}\label{eq:blocked}
 \begin{pmatrix}
 -\frac{1}{2m} \Delta_{l_1,l_1}  + V_{l_1m_1\,l_1m_1;l_1m_1 l_1m_1}-E &  V_{l_1m_1\,l_1m1;l_1m1l_2m_2} & \ldots \\
       V_{l_1m_1\,l_2m_2;l_1m_1\,l_1m_1} & -\frac{1}{2m}\Delta_{l_1,l_2} + V_{l_1m_1\,l_2m_1;l_1m_1l_2m_2} -E &\ldots\\
\vdots & \vdots & \ddots
 \end{pmatrix},
\end{equation}}
where the $\Delta_{l_1,l_2}$ are the radial differential operators
defined in \eqref{eq:definition_delta}. This can be written as a
coupled Helmholtz operator  \wv{\begin{equation}
 \begin{pmatrix}
 -\Delta_{l_1,l_1}  -k^2_{l_1m_1\,l_1m_1;l_1m_1 l_1m_1}(\rho_1,\rho_2)&  -k^2_{l_1m_1\,l_1m1;l_1m1l_2m_2}(\rho_1,\rho_2) & \ldots \\
       -k^2_{l_1m_1\,l_2m_2;l_1m_1\,l_1m_1}(\rho_1,\rho_2) & -\Delta_{l_1,l_2} - k^2_{l_1m_1\,l_2m_1;l_1m_1l_2m_2}(\rho_1,\rho_2)&\ldots\\
\vdots & \vdots & \ddots
 \end{pmatrix},
\label{eq:blocked_helmholtz}
\end{equation}}

After discretization of the differential operators on a grid \hbz{discussed in detail in Section \ref{sec:discretize}, we arrive at a system of linear equations,} $Ax=b$. The matrix \hbz{$A$} will have the same blocked structure as the coupled system of partial differential equations above and we can write:
\begin{equation}
 \begin{pmatrix}
  A_{11}   & A_{12}  &  A_{13} & \ldots \\
  A_{21}   & A_{22}  &    &  \\
  A_{21}   &      &  A_{33} &  \\
  \hdots  &    & & \ddots
 \end{pmatrix}
\begin{pmatrix}
x_1\\
x_2 \\
x_3\\
\vdots
\end{pmatrix}
= \begin{pmatrix}
  b_1\\
  b_2 \\
  b_3\\
\vdots
  \end{pmatrix},
\end{equation}
where the discretized differential operators will only appear in the
diagonal blocks $A_{ii}$. Since the differential operators
will lead to the largest eigenvalues, the condition number of the full matrix $A$ will also be determined
predominantly by the diagonal blocks $A_{ii}$. \edt{
After discretization of
$\rho_1$ on $n$ grid points and $\rho_2$ on $n$ grid points, a single
block is a sparse matrix of size $n^2 \times n^2$. }

\hbz{The solution method for solving this type of breakup problems as employed in \cite{Wim05} is iterative. This method was developed in \cite{baertschy2001solution} and exploits the particular block structure of $A$.} A block diagonal preconditioning matrix $M$ is constructed that contains only the diagonal blocks $A_{ii}$.  Since the largest eigenvalues and eigenvectors of $M$ and $A$ are very similar, $M^{-1}A$ has a smaller condition number, \hbz{and therefore, $M^{-1}$ proves to be a good preconditioner for any suitable Krylov-subspace method. However, note here in particular, that the strategy in used in \cite{Wim05} is to exactly invert the blocks \edt{(each of size $n^2 \times n^2$)} within the preconditioning step. Inasmuch as each diagonal block represents a two-dimensional system, the diagonal block matrices can be inverted \edt{possibly} on a single processor. The coupled system, however, requires the inversion of many diagonal blocks and requires a cluster.}

\wv{
However, the problems currently under investigation in the physics and
chemistry communities such as the impact-ionization problems or
problems where electronic motion and nuclear motion are combined
described in section \ref{sec:impactionisation}, each diagonal block
\hbz{constitutes} a three dimensional problem and \hbz{renders itself too unwieldy for exact inversion}}. 

We therefore study in this paper the \hbz{multigrid-preconditioned iterative solution of the diagonal block only, and not the entire problem as a whole. The diagonal blocks \eqref{eq:blocked} corresponds closely to the model problem \eqref{eq:modelproblem}. It is important to understand that the complete process now involves two independent iterative schemes, the outer scheme for approximately inverting the entire system via a preconditioned Krylov process, and the inner iterative scheme which uses multigrid preconditioning for approximately inverting the diagonal blocks within the outer preconditioner. The latter alone forms the subject matter of this paper. We have chosen to restrict the dimensions to two for this study.}

\subsection{Examples}
\wv{To illustrate the significance of the coupled system of partial
differential equations we give a few example physical systems that are
currently studied with the approach. We cite the relevant papers.}
\subsubsection{The dynamics of two electrons in  a Helium atom}
The Helium atom is a quantum system that has two electrons with a negative charge
and one nucleus \hbz{which} has a positive charge of unit two.
Since the nucleus is much heavier than each electron, the position of the
nucleus is taken as the center of the coordinate system.
In this coordinate system the first electron is at $\mathbold{r}_1$ and the second electron at
$\mathbold{r}_2$ The potentials in equation \eqref{eq:schrodinger} are then
\begin{equation}
\begin{aligned}
  V_1(\mathbold{r}_1) = -\frac{2}{|\mathbold{r}_1|},\quad  V_2(\mathbold{r}_2) = -\frac{2}{|\mathbold{r}_2|} ,\quad   V_{12}(\mathbold{r}_1,\mathbold{r}_2) = \frac{1}{|\mathbold{r}_1-\mathbold{r}_2|}
\end{aligned}
\end{equation}
To arrive at the potentials in the coupled problem \eqref{eq:coupled_potential} the multipole expansion
\begin{equation}
 \frac{1}{|\mathbold{r}_1-\mathbold{r}_2|} = \sum_l \frac{\rho_<^l}{\rho_>^{l+1}} P_l(\cos(\theta_{12}))  \label{eq:multipole}
\end{equation}
is used to expand $V_{12}$. Where $\rho_<$ and $\rho_>$ denote,
respectively, the smallest and largest of $\rho_1$ and $\rho_2$. The
angle $\theta_{12}$ is between the vectors $\mathbold{r}_1$ and
$\mathbold{r}_2$.  Since $V_1$ and $V_2$ are central potentials they
will appear as $-2/\rho_1$ and $-2/\rho_2$ on the diagonal blocks of
\eqref{eq:coupled_potential} when $l_1=l_1^\prime$ and
$l_1=l_1^\prime$.  The multipole expansion \eqref{eq:multipole},
however, will lead to potentials that couple the blocks with different
$l$ values in equations \eqref{eq:coupled}.  Since the problem is
symmetric around the $z$ axis, different $m$ blocks are
decoupled. Recent processes in Helium studied with this approach are
one and two-photon double ionization
\cite{mccurdy2004theoretical,horner2007two}.

\subsubsection{The dynamics of two electrons in the Hydrogen molecule}
The Hydrogen molecule consist of two negatively charged electrons and
two protons with a positive charge. The two protons are much heavier
than the electrons.  After the Born-Oppenheimer approximation \wv{the two protons can be considered fixed in space}. The
dynamics of the two electrons are governed by equation
\eqref{eq:schrodinger}, where the potential is given by the static
field of the charged protons.  If we take a coordinate system around
the center of mass of the protons and $\mathbold{R}$ is the vector
connecting the two protons \wv{and $\mathbf{r}_1$ and $\mathbf{r}_2$ the coordinates of the electrons}, the potentials in \eqref{eq:schrodinger}
are
\begin{equation}
 \begin{aligned}
  V_1(\mathbold{r}_1) &= -\frac{1}{|\mathbold{r}_1-\mathbold{R}/2|} -\frac{1}{|\mathbold{r}_1+\mathbold{R}/2|}, \quad
  V_2(\mathbold{r}_2) = -\frac{1}{|\mathbold{r}_2-\mathbold{R}/2|} -\frac{1}{|\mathbold{r}_2+\mathbold{R}/2|},
\end{aligned}
\end{equation}
\begin{equation}
 \begin{aligned}
 V_{12}(\mathbold{r}_1,\mathbold{r}_2) &= \frac{1}{|\mathbold{r}_1-\mathbold{r}_2|}
 \end{aligned}
\end{equation}
The first is the attraction of the first electron to the two
protons. The second is the same attraction but for the second
electron.  The third potential is the electron-electron repulsion
because both have a negative charge.  To derive the potentials in the
coupled basis we use, again, the multipole expansion
\eqref{eq:multipole} for each of the potentials.  \wv{Now all potentials 
couple the blocks.}  Again, an example of a process studied in this
approach is one-photon double ionization \cite{Wim05}.

\subsubsection{Impact ionization of Helium or the Hydrogen molecule}\label{sec:impactionisation}
  When an \wv{additional} electron with sufficient energy collides and breaks up the
  Helium atom or Hydrogen molecule \wv{that have already two electrons, we are tracking three particles}. We then have a 9D problem.
  If we denote the coordinate of the third, impacting, electron as
  $\mathbold{r}_3$ we end up with Helmholtz operator
\begin{equation}
  -\Delta_{\mathbold{r}_1} - \Delta_{\mathbold{r}_2} - \Delta_{\mathbold{r}_3} - k^2(\mathbold{r}_1,\mathbold{r}_2,\mathbold{r}_3).
\end{equation}
After \wv{the partial wave expansion} we end up, again, with a coupled problem 
\begin{equation}
\begin{aligned}
&\Bigg[ -\frac{1}{2m}\frac{\partial^2}{\partial \rho_1^2} -\frac{1}{2m}\frac{\partial^2}{\partial \rho_2^2} 
 -\frac{1}{2m}\frac{\partial^2}{\partial \rho_3^2} + \frac{l_1(l_1+1)}{2 \rho_1^2}  + \frac{l_2(l_2+1)}{2 \rho_2^2}  + \frac{l_3(l_3+1)}{2 \rho_3^2} \\
   &+V(\rho_1,\rho_2,\rho_3)-E \Bigg]  \psi(\rho_1,\rho_2, \rho_3) 
 + \sum_{l_1 \neq l_1^\prime,l_2\neq l_2^\prime,l_3\neq l_3^\prime} \!\!\!\!\!C_{l_1l_2l_3,l^\prime_1l^\prime_2l^\prime_3}
(\rho_1,\rho_2,\rho_3) 
\psi_{l_1^\prime,l_2^\prime,l_3^\prime}(\rho_1,\rho_2,\rho_3)  = 0
\end{aligned}\label{eq:coupled3electron}
\end{equation}
where $V(\rho_1,\rho_2,\rho_3)$ is a diagonal block potential while
$C(\rho_1,\rho_2,\rho_3)$ couples the blocks.  Because of the scale of
these problem there are currently no converged results for this
problem.  A similar high dimensional problem can be formulated for
Hydrogen molecule when no Born-Oppenheimer approximation is applied.
Then the motion of the electrons, $\mathbold{r}_1$ and
$\mathbold{r}_2$, is coupled to the motion of the protons
$\mathbold{R}$.  \wv{Solving these problems is a great interest in the
scientific community.}

\section{Discretization}\label{sec:discretize}
\subsection{Absorbing boundary conditions}
In order to solve equations, such as the Helmholtz equation, defined on an unbounded domain $\Omega_0\subseteq\mathbb{R}^d$ numerically, an assumption is made on the asymptotic behavior of the solution. The truncated computational domain is a bounded subset $\Omega\subset\Omega_0$ of the original one, with artificially introduced boundary conditions that imply the postulated asymptotic behavior. A commonly used example are the first order Sommerfeld radiation boundary conditions applied to the homogeneous Helmholtz problem, $\frac{\partial u}{\partial \hat{n}} + i k u = 0$, where $\hat{n}$ is the outward normal on the boundary $\partial\Omega$. An exponential decay of the solution depending on the constant wave number $k$ is assumed towards the boundary.

In more complicated Helmholtz models such as the one derived from the Schr\"odinger equation \eqref{eq:modelproblem} more robust techniques are preferable. In the \emph{perfectly matched layer} (PML) approach \cite{B94} a small boundary layer $\Gamma\subset\mathbb{R}^d$ is added beyond any point of domain truncation. On this finite layer the continuous model is adapted to capture the asymptotic behavior, with trivial boundary conditions at the end of the layers $\partial\Gamma$. This idea is equivalent to a complex coordinate stretching \cite{CW94,R95,TC98} in the boundary layers, where the original equation is used in the new coordinate system $\Gamma_z\subset\mathbb{C}^d$ with homogeneous Dirichlet boundary conditions at the end $\partial\Gamma_z$, also known as \emph{exterior complex scaling} (ECS) \cite{NB78,S79}. In general we can define an analytic continuation on the layers by
\begin{align*}
 z(x) = \left\{
  \begin{array}{ll}
    x, & \hbox{$x \in \Omega$;} \\
    x+if(x), & \hbox{$x \in \Gamma$,}
  \end{array}
\right.
\end{align*}
with $f \in \mathcal{C}^2$ increasing (e.g. linear, quadratic, \ldots) and $\displaystyle \lim_{x\to\partial\Omega}f(x) = 0$. We
denote the image of the layer $\Gamma_z\equiv z(\Gamma)$ and call it the \emph{complex contour}. This boundary layer method does not need an explicit input of the wave number and it can easily be tuned in numerical experiments. Because of the straightforward mathematical meaning the ECS method is interesting in numerical analysis.

\subsection{Finite difference }
ECS boundary conditions and their application in chemical reactions
have been used in finite difference, B-spline and spectral element
discretization \cite{mccurdyTR2004}. \wv{Finite Elments methods are
  hardly used for this type of problems because the computational
  domain is often a square or a rectangular strip.} In this article we
use finite differences \wv{since this low order discretisation can
  already help us to understand the convergence of the iterative
  method.} We define a one-dimensional uniform grid $(z_j)_{0\leq
  j\leq n}$ on the real interval $[0,1]$ with $z_0 = 0$ and $z_n=1$
and mesh width $h = 1/n \in\mathbb{R}$. Starting in $1$, we apply
linear ECS, so the absorbing layer is a line connecting $1$ and
$R_z\in\mathbb{C}$ henceforth denoted by $[1,R_z]$. A second uniform
grid $(z_j)_{n\leq j\leq n+m}$ discretized this complex contour, with
$z_{n+m} = R_z$ and complex mesh width $h_\gamma = (R_z-1)/m$. The
union of these two grids is the ECS grid
\begin{align}\label{eq:ecsgrid}
(z_j)_{0\leq j\leq n+m} \mbox{ on } [0,1]\cup[1,R_z]
\end{align}
in the entire ECS domain. We will denote the fraction $\gamma = h_\gamma/h = (R_z-1)/(R-1)$.
\begin{figure}[htbp!]
\begin{center}
\includegraphics[width=12cm]{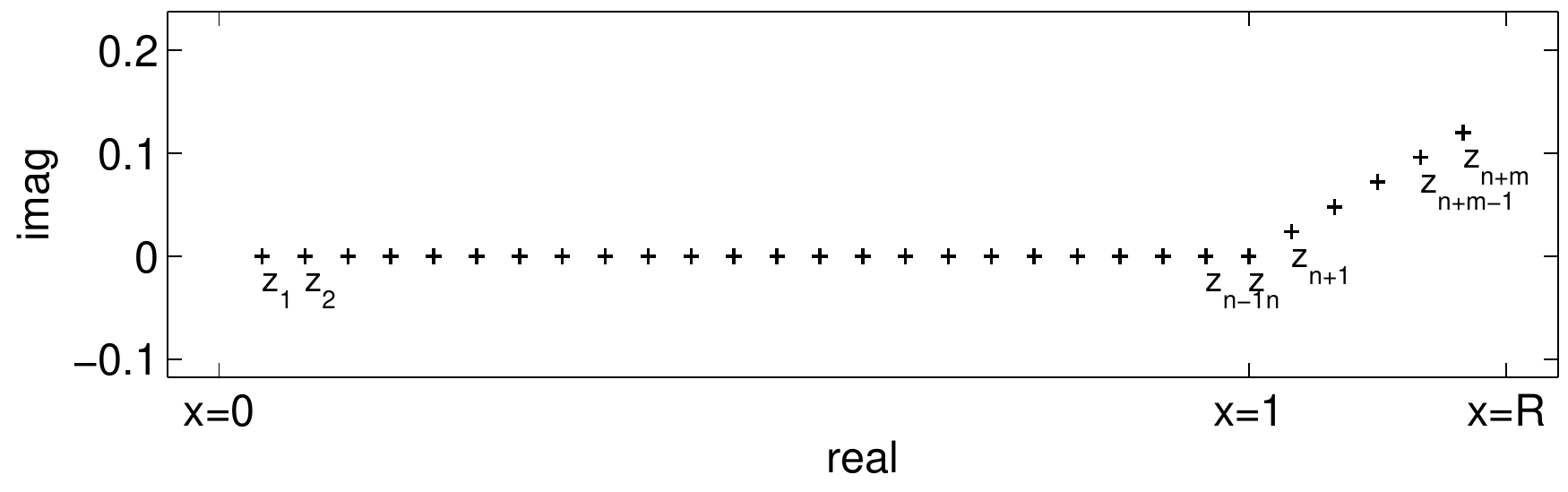}\caption{Discretized ECS domain $z_j$. The ECS domain is discretized with complex mesh widths on the complex contour.}\label{fig:disecs}
\end{center}
\end{figure}

A thorough numerical analysis of the negative Laplace operator $L=-\Delta$ discretized on this ECS domain yields some important insights for the use of ECS on more general operators. To approximate the second derivative we \hbz{employ the following standard formula for un-equal mesh sizes, and non-uniform grids:}
\begin{align*}
\frac{d^2u}{d z^2}(z_j) \approx
\frac{2}{h_{j-1}+h_j}\left(\frac{1}{h_{j-1}}u_{j-1}-\left(\frac{1}{h_{j-1}}+\frac{1}{h_j}\right)u_j
+\frac{1}{h_j}u_{j+1}\right)
\end{align*}
for non-uniform grids in grid point $j$, where $h_{j-1}$ and $h_j$ are the left and right mesh widths respectively, and may belong either to the $h$ category or to the $h_\gamma$ category. \wv{The formula can be easily derived as \edt{an} exercise in Taylor expansions} and it reduces to regular second order central differences when $h_{j-1}=h_j$, i.e., in the interior real region $(0,1)$, and in the interior of the complex contour $(1,R_z)$ because the scaling function $f$ is taken to be linear. The only exception is the point $z_n$ where at most we lose an order of accuracy, however with ample discretization steps, the overall accuracy is anticipated to match up to second order. We will denote the resulting discretization matrix $L_h$.

\subsection{Spectral properties}
\wv{The hardest model problem, from a\edt{n iterative point} of view, is 
the problem with $l_1=0$ and $l_2=0$ since the problem is then \edt{at its} most
indefinite \edt{state}. For larger $l_1$ and $l_2$ the problem becomes more
definite. Therefore we focus on the remainder of the paper on the
problem with $l_1=0$ and $l_2=0$.}

The spectrum of the discretization matrix $L_h$ determines the
convergence behavior of iterative methods such as Krylov subspace
methods and multigrid schemes for solving any system $L_h u_h =
b_h$. The spectrum $\sigma(L_h)$ is drastically different from the
spectrum $\sigma(L)$ of the continuous operator, on the undiscretized
ECS grid $[0,1]\cup[1,R_z]$. Indeed, $\sigma(L) = \left\{\left(\frac{j
  \pi}{R_z}\right)^2 |j\in\mathbb{N}_0 \right\}$, is an infinite set
of points on the complex line $\rho e^{-i2\theta_\alpha}$, with
$\rho\in\mathbb{R^+}$ and $\theta_\alpha$ the complex angle of the
complex boundary $R_z$. The shape of the spectrum $\sigma(L_h)$ of the
discretization matrix is less obvious as follows from the next lemma,
that is proved in \cite{reps2009}.

\begin{lem}\label{lem:eigdis}
Consider the ECS grid \eqref{eq:ecsgrid} and the discretization matrix $L_h$. Define $\gamma =
\frac{h_\gamma}{h}$. Then the eigenvalues of $L_h$ are the solutions of
\begin{align*}
F(\lambda) \equiv \frac{\tan(2n p(\lambda))}{\tan(2m q(\lambda))}+\frac{\cos(p(\lambda))}{\cos(q(\lambda))} = 0,
\end{align*}
with $p(\lambda)=\frac{1}{2}\arccos(1-\frac{\lambda}{2}h^2)$, $q(\lambda)=\frac{1}{2}\arccos(1-\frac{\lambda}{2}\gamma^2
h^2)$.
\end{lem}

\begin{figure}[h!]
\begin{center}
\includegraphics[width=9cm]{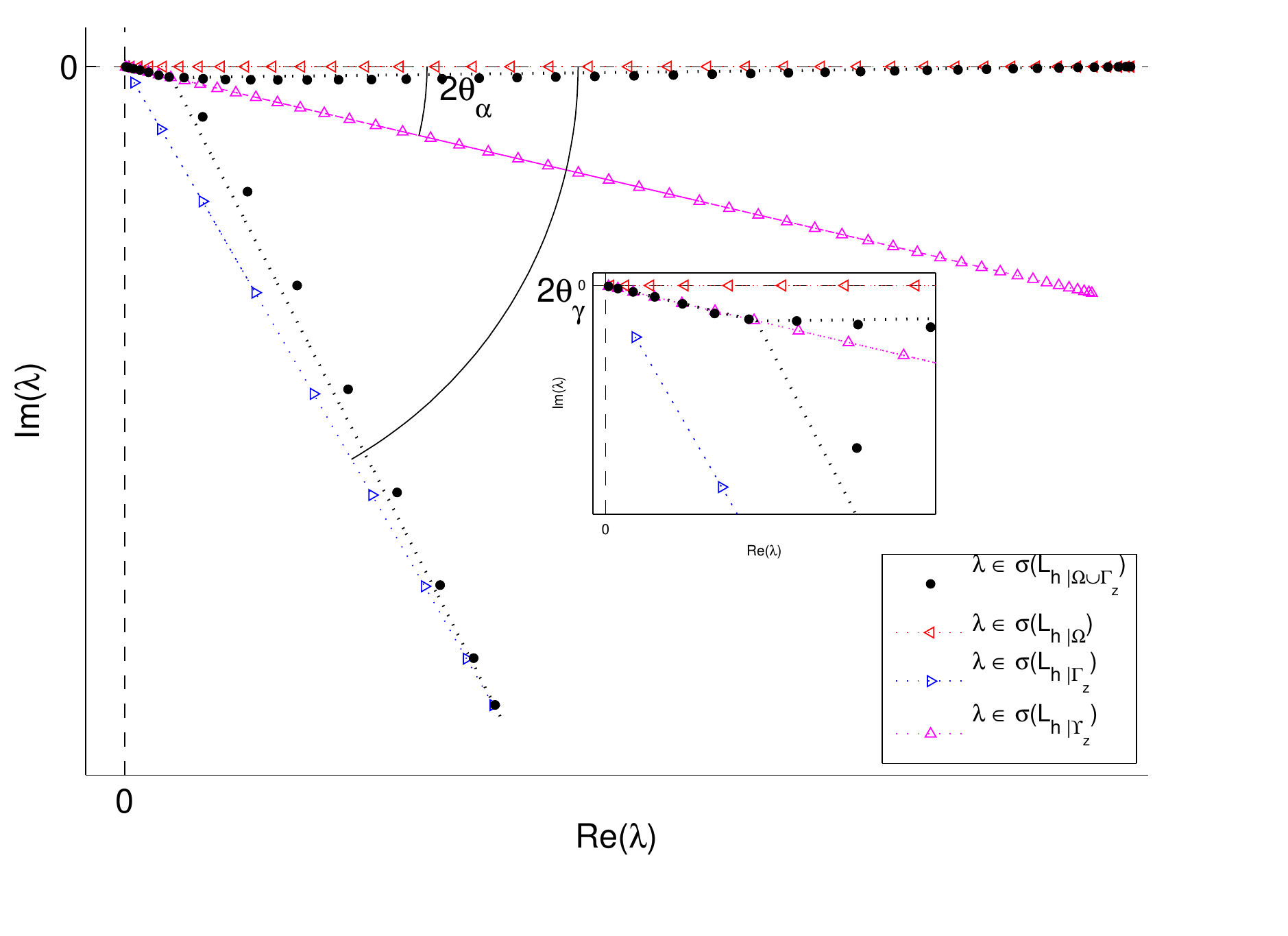}\caption{
\wv{ The eigenvalues of the ECS Laplacian discretization matrix ($\bullet$) lie along a pitchfork with a $Y$ shape, a result in Lemma~\ref{lem:eigdis}. 
 A part of the eigenvalues lie close to the eigenvalues the same Laplace problem restricted to the interior real domain ($\triangleleft$).  In a similar way part of the eigenvalues lie close to the eigenvalues when the Laplacian is restricted to the complex contour ($\triangleright$).} \wv{The inset shows the area around the origin where the smallest eigenvalues are approximated by the smallest eigenvalues of the Laplace problem defined on the complex line }$[0,R_z]$ ($\triangle$).
 (color online)}
\label{fig:dowsingrod}
\end{center}
\end{figure}

For the Laplace problem the ECS discretized spectrum has the typical Y-shape of a pitchfork. There is a clear complex branch associated to eigenvectors located on the complex contour, along the complex line $[0,4/h_\gamma^2]$, and a branch closer to the real line $[0,4/h^2]$ that corresponds to eigenvectors located on the real domain. The smallest eigenvalues, in the small tail of the pitchfork, belong to the smoothest eigenvectors spread over the entire ECS domain. They lie close to the smallest \edt{eigenvalue} of the continuous ECS operator $L$ (Figure \ref{fig:dowsingrod}), that is along the complex line $[0,4/h_\alpha^2]$, where $h_\alpha = R_z/(n+m+1)$ is the complex mesh width, belonging to a straight complex grid connecting $0$ and $R_z$.

\section{The $QD$ preconditioner and multigrid} \label{sec:precmg}
\subsection{The preconditioner} \label{subsec:precon}
We use a preconditioning operator that has a spectrum bounded by a single quadrant such that it 
can be approximately inverted with standard multigrid components, which are clearly unstable for indefinite problems.

In this article we compare the use of a preconditioner $\tilde{\mathcal{Z}}$ which is a scaled 
and shifted version of the original Helmholtz operator $\mathcal{Z} \wv{ = -\Delta_{l_1,l_2}-k^2(x,y)}$ defined in equation \eqref{eq:modelproblem}.
We propose to use
\begin{equation}
\tilde{ \mathcal{Z}} = \delta^2 \mathcal{Z}  + (1-i\delta k)
  \label{eq:preconditioner_operator}
\end{equation}
where $\delta \in \mathbb{R}$ is chosen such that $\tilde{\mathcal{Z}}$ is definite. This preconditioner is very similar to the
one proposed in \cite{meerbergen2009}.
 Suppose $\lambda_0$ is that eigenvalue of the original operator $\mathcal{Z}$ which has the smallest real part.
We can then choose $\delta$ such that $-\delta |\lambda_0| + 1 \ge 0$.  For a Helmholtz problem with a constant wave number $k$
\wv{and $l_1=0$ and $l_2=0$} this would mean that $\delta \le 1/k$.



The eigenvalues of the preconditioned operator \wv{$\tilde{\mathcal{Z}}^{-1} \mathcal{Z}$} lie inside a circle of radius
$1/(2\delta^2)\left|1-i/(k\delta)\right|$ centered around
$\frac{1}{\delta^2}\left(\frac{1}{2}-i\frac{1}{2k\delta}\right)$.  We
can readily see this with the following arguments. The preconditioner
$\tilde{\mathcal{Z}}=\delta^2\mathcal{Z}+(1-i\delta k)$ has the same
eigenvectors as the $\mathcal{Z}$. The eigenvalues of the
preconditioned system \wv{$\tilde{\mathcal{Z}}^{-1}\mathcal{Z}$} are therefore given by
\begin{align}\label{eq:eigMA}
\frac{\lambda}{\delta^2\lambda+(1-i\delta k)}, 
\end{align}
where $\lambda$ is an eigenvalue of \wv{$\mathcal{Z}$}. We assume that
the eigenvalues of $\mathcal{Z}$ are located in the lower half of the complex
plane, $\sigma(\mathcal{Z})\subset\mathbb{C}^-$. Then $\sigma(\tilde{\mathcal{Z}}^{-1}\mathcal{Z})$ is
inside the circle that is the image of the real axis of the transform
\eqref{eq:eigMA}. This circle goes through $0$,
$1/\delta^2\in\mathbb{R}$ and $-i/(k\delta^3)\in i\mathbb{R}$, so the
center $c$ is the crossing point of the lines $\Re(z)=1/(2\delta^2)$
and $\Im(z)=-1/(2k\delta^3)$,
\begin{align*}
c = \frac{1}{\delta^2}\left(\frac{1}{2}-i\frac{1}{2k\delta}\right).
\end{align*}
And so the radius is $r=|c|= 1/(2\delta^2)\left|1-i/(k\delta)\right|$.

\subsection{Multigrid}
\label{subsec:mg}
Heuristically, we note that multigrid (for convergence and efficiency), has more stringent requirements on 
the condition number of the spectrum when it crosses into different quadrants of the complex plane, than when 
it does not. The preconditioner in Section \ref{subsec:precon} has the property that its spectrum is 
restricted to the fourth quadrant. It can therefore be very efficiently inverted by multigrid using the standard 
components, which include $\omega$-Red Black Jacobi, with $\omega = 1.05$, Full Weighting averaging for restriction, 
Bilinear interpolation for prolongation, and rediscretization on the coarse grids; in a simple V(1,1) cycle set-up.
 For experimental purposes we compare the performance of this quadrant-definite ($QD$) preconditioner with the $CSG$ and 
the $CSL$ preconditioners.
 In \cite{reps2009}, we showed that the $CSL$ and the $CSG$ preconditioner can be inverted efficiently using multigrid based on matrix components, such as ILU-smoother and the Galerkin coarse grid operator in a V(0,1) cycle set-up. This study is more focused on using matrix-free components, which leaves only the $\omega$-Jacobi smoother, and the discretization coarse grid operator for the $CSL$ and the $CSG$ preconditioners. Moreover, this multigrid has to be employed in an F$_{\gamma_f}^{\gamma_c}$(1,1) cycle set-up. 

\begin{algorithm}
\centering {$u_l^{m+1} = {\rm \textbf{MG}}(l, \bsA_l, b_l, u_{l}^m, \nu_1, \nu_2, C, \gamma_f, \gamma_c)$.}
\caption{\bf ~~ Multigrid pseudocode}
\label{alg:mg}
\begin{enumerate} \sf 
\item[\bf (0)] \sf  {\bf Initialization} \\[1.0ex]
--  \sf  \underline{If} $l=C$, $u_l^{m+1} = $ \textbf{exact} ($\bsA_l, b_l$); \underline{Bail out}; \underline{endif}\\
--  \sf  Build the coarse-grid operator $\bsA_{l+1}$, and the restriction $I_l^{l+1}$, and prolongation $I_{l+1}^l$ operators.\\
\item[\bf (1)] {\bf Pre--smoothing} 
\begin{itemize}
\item[--] \sf  Compute
$u^{m+\nicefrac{1}{3}}_l$  by applying $\nu_1 ( \ge 0)$ smoothing
steps to $u^m_l:$\\
$  
u^{m+\nicefrac{1}{3}}_l = {\rm \textbf{smooth}}^{\normalsize\mbox{$\nu_1$}} ( \bsA_l,
b_l, u^m_l)\enspace.
$
\end{itemize}
\item[\bf (2)] \sf  {\bf Coarse grid correction} \\[1.0ex]
\begin{tabular}{lll}
-- &   Compute the residual  & $  \ol{r}^m_l = b_l - \bsA_l u^{m+\nicefrac{1}{3}}_l$ \enspace.\\[1.0ex]
-- &   Restrict the residual & $  \ol{r}^m_{l+1} = I_l^{l+1} \hspace{2mm} \ol{r}^m_l$ \enspace.\\[1.0ex]
-- &   Compute the approximate error \\
   &   $\widehat{e}^m_{l+1}$ from the \textit{defect equation}. & $\bsA_{l+1} \hspace{2mm} \widehat{e}^m_{l+1} = \ol{r}^m_{l+1}$\\
   &   & \\	
   &	by the following recursion& \\
\end{tabular}

\hspace*{5ex}\framebox[10cm]{
\begin{minipage}{8cm}
\vspace{2mm}
 \sf 
$\text{$\hspace{-7mm}$}$\underline{If} $l=C$, $\hat{e}_{l+1}^{m} = $ \textbf{exact} ($\bsA_{l+1}, \ol{r}_{l+1}^m$);  \underline{endif}\\ \sf 
$\text{$\hspace{-7mm}$}$\underline{If} $l<C$, approximate $\hat{e}_{l+1}^{m}$ recursively:\\
$\text{$\hspace{-5mm}$}$ $\hat{e}_{l+1}^{m,1} = \ol{0}$; \\
$\text{$\hspace{-4mm}$}$\underline{do} $i = 1$ to $\gamma_c$  \\
$\text{$\hspace{-2mm}$}$ \underline{If} $\hspace{1mm}$ $i == 1$,\\
$\text{$\hspace{1mm}$}$ $\hat{e}_{l+1}^{m,i+1}={\rm\textbf{MG}}(l+1,\bsA_{l+1},\ol{r}_{l+1}^m, \hat{e}_{l+1}^{m,i},\nu_1,\nu_2, C, \gamma_f, \gamma_c)$\\
$\text{$\hspace{-2mm}$}$ \underline{else} \\
$\text{$\hspace{1mm}$}$ $\hat{e}_{l+1}^{m,i+1}={\rm\textbf{MG}}(l+1,\bsA_{l+1},\ol{r}_{l+1}^m, \hat{e}_{l+1}^{m,i},\nu_1,\nu_2, C, \gamma_c, \gamma_f)$\\
$\text{$\hspace{-2mm}$}$ \underline{endif}\\
$\text{$\hspace{-4mm}$}$\underline{continue} i\\
$\text{$\hspace{-7mm}$}$\underline{endif}
\end{minipage}}
\vspace*{2mm}

\begin{tabular}{lll}
-- & \sf  Interpolate the correction \hspace*{8ex} 
      &$  \widehat{e}^m_l = I_{l+1}^l \hspace{2mm} \widehat{e}^m_{l+1}$ \enspace.\\[1.0ex]
-- & \sf  Compute the corrected \\
   & \sf  approximation on $\Omega_l$ 
                & $  u_l^{m+\nicefrac{2}{3}} = u^{m+\nicefrac{1}{3}}_l + \widehat{e}^m_l$\enspace.
\end{tabular}
\\
\item[\bf (3)]  \sf 
{\bf Post--smoothing}
\begin{itemize}
\item[--] \sf  Compute $u^{m+1}_l$ by applying $\nu_2\; (\ge \! 0)$ smoothing
        steps to $u_l^{m+\nicefrac{2}{3}}$:\\
$   u^{m+1}_l = {\rm \textbf{smooth}} ^{\normalsize\mbox{$\nu_2$}} 
    ( u_l^{m+\nicefrac{1}{2}} , \bsA_l, b_l) \enspace.
$
\end{itemize}
\end{enumerate} 
\end{algorithm}

\begin{figure}
\centering
 \psfig{figure=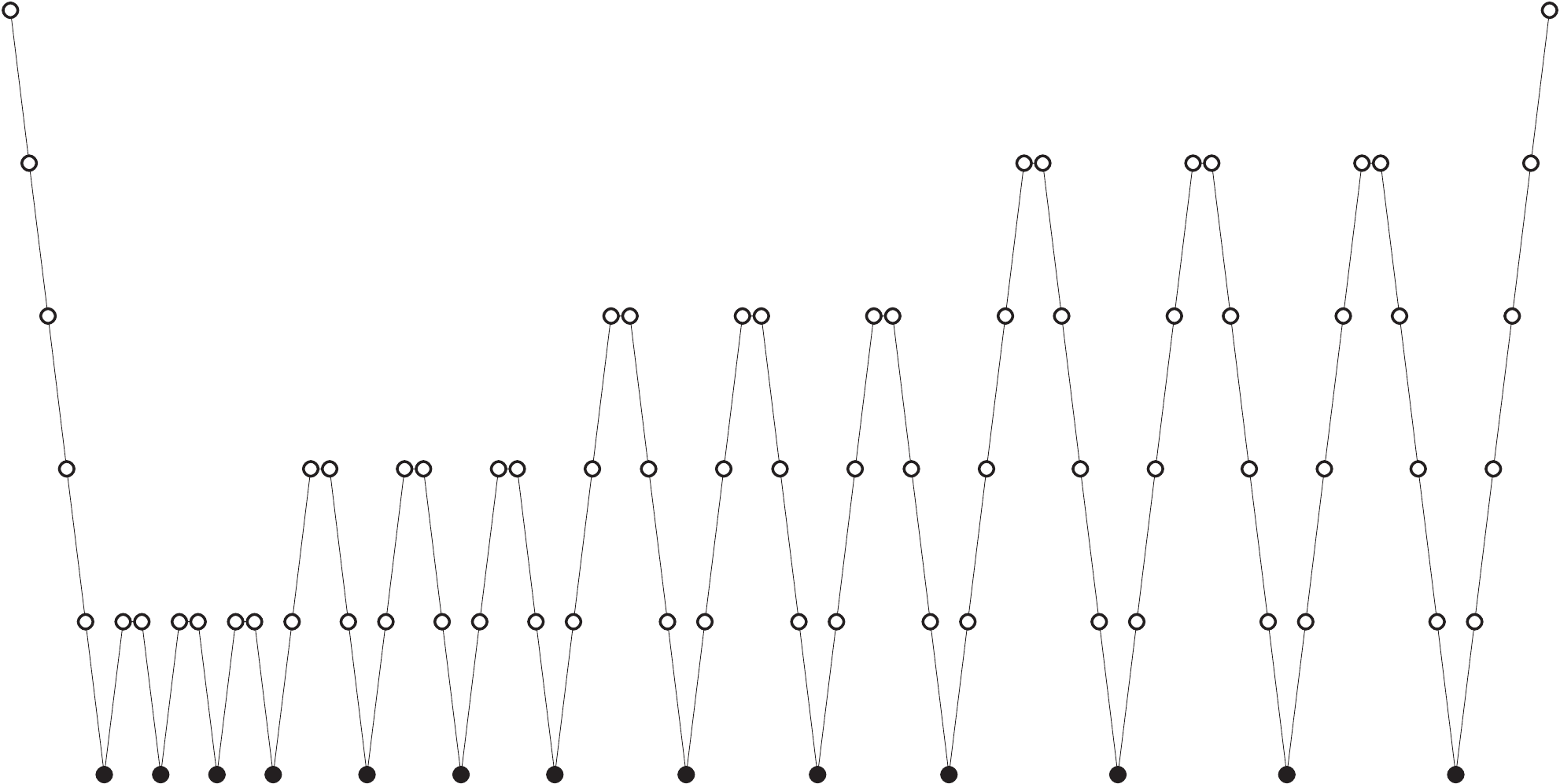, width=10cm}
\caption{An F$_1^4$(1,1) cycle, obtained by Algorithm \ref{alg:mg}, with
 $\gamma_f=1$ and $\gamma_c=4$. $\circ$ stands for smoothing, $\backslash$ stands for
 restriction to a lower level, / stands for prolongation to a higher level, $\bullet$ stands for exact solution.}
\label{fig:F14_cyc}
\end{figure}

Algorithm \ref{alg:mg} is a slight variation of the standard multigrid algorithm in \cite{Trot01}, 
and yields various cycle types depending on the values of $\gamma_f$ and $\gamma_c$ in a unified manner. 
In this algorithm, $l$ indicates the current level, $C$ the coarsest level, and $\bsA$, the discrete operator 
(i.e., $CSL$/$CSG$/$QD$ precond. operators) at various levels. $b_l$ is the right hand side, $u_l^m$ is the starting
 guess, and $\nu_1$ and $\nu_2$ are the number of pre and post smoothing sweeps. $\gamma_f$ and $\gamma_c$ are
 the cycle indices used at the fine and the coarse levels, respectively. E.g., calling this method with 
$\gamma_f=1$ and $\gamma_c=1$ gives the standard V cycle, $\gamma_f=2$ and $\gamma_c=2$ gives the standard W cycle, 
while $\gamma_f=1$ and $\gamma_c=2$ yields the standard F cycle. $\gamma_f=1$, and $\gamma_c=n$ renders an F cycle 
with $n-1$ recursions on the coarse levels. We found that in the context of inverting $CSL$ and $CSG$ based Helmholtz
preconditioners, F cycles with 2 and 3 recursions on the coarse levels were particularly beneficial over the 
standard F cycle. F cycle with 3 recursions is abbreviated as F$_1^4$ and is shown in Figure \ref{fig:F14_cyc}.

\section{Numerical Results}
\label{sec:numex}
 In this section we conduct numerical experiments on the Helmholtz
 model problems given by the generic prototype
 $\mathcal{Z}u(x,y)=f(x,y)$. \wv{Again we focus on problems with $l_1=0$
 and $l_2=0$ because these are the hardest problems.} The operator
 $\mathcal{Z}$ is defined in each of the following case as follows:
\begin{align}
 \mathcal{Z} &= -\Delta - k^2 & \text{MP1} \\
 \mathcal{Z} &= -\Delta - \nu(\frac{1}{e^{x^2}}+\frac{1}{e^{y^2}}) - k^2 & \text{MP2}\\
 \mathcal{Z} &= -\Delta -\frac{1}{x} - \frac{1}{y} - k^2 & \text{MP3}
\end{align}
 
For MP1, $f(x,y)$ is the Dirac delta function that stays zero throughout the domain save one point in the middle where it assumes the value 1. For MP2 and MP3, $f(x,y) = e^{x^2+y^2}$. We use ECS-ABL on all four sides with MP1, and on the north and the east sides only with MP2 and MP3. \hbz{The model problems are solved iteratively with Bi-CGSTAB preconditioned with multigrid approximated inverses of the following operators}:
\begin{align}
 \mathcal{M}_{CSL} &= -\Delta + (\beta_1 + i\beta_2)k^2,\\
 \mathcal{M}_{CSG} &= Z, & {\scriptstyle \text{on the grid rotated by angle $\theta_{\alpha}$ in the complex plane, see \cite{reps2009}}}\\
 \mathcal{M}_{QD} &= (1-i)I + \frac{1}{\lvert \text{Re}(\lambda_0) \rvert}Z & {\scriptstyle \text{$\lambda_0$ is the eigenvalue of $Z$ with the smallest real part.}}
\end{align}

\hbz{$M_{CSL}$ is the \textit{Complex-shifted Laplacian} as appears in \cite{Yogi06}. A small complex shift is added to the Laplacian operator. This imparts a rectangular translational effect on the operator spectrum. $M_{CSG}$ is the original operator discretized on the so-called \textit{Complex-scaled Grid} and appears in detail in \cite{reps2009}. The basic mesh size has been multiplied with $e^{i\theta_{\alpha}}$. This imparts a rotation to the operator spectrum about the origin by an angle equal to $\theta_{\alpha}$. $M_{CSG}$ is as efficient as $M_{CSL}$ in general, and slightly better for the current problems. }

\begin{table}
 \begin{tabular}{|c|c|c|c|c|} \hline 
  \multirow{2}{*}{Preconditioner} & \multicolumn{2}{c|}{Multigrid} & mg cyc. & Bi-CGSTAB\\ \cline{2-3}
 & cyc, smooth., $\omega$ & mg-conv. , \# cycles & per prec. & iter , cputime \\ \hline
$CSL$ & F$_1^4$(1,1), $\omega$-Jacobi & 0.43 , 17 & \multirow{2}{*}{1} & \multirow{2}{*}{60 , 2m 11s}\\
$(\beta_1, \beta_2 = (-1, -0.3)$ & $\omega = 0.8$ & 4.21s & & \\ \hline
$CSG$ & F$_1^3$(1,1), $\omega$-Jacobi & 0.39 , 15 & \multirow{2}{*}{1} & \multirow{2}{*}{62 , 2m 2s}\\
$\theta_{\alpha} = \nicefrac{\pi}{14}$ & $\omega = 0.8$ & 3.18s & & \\ \hline
$QD$ & V(1,1), $\omega$-RB Jacobi & 0.09 , 6 & \multirow{2}{*}{1} & \multirow{2}{*}{170 , 5m 39s}\\
$\text{Re}(\lambda_0) = -2.6 \times 10^4$ & $\omega = 1.0$ & 1.2s & & \\ \hline
 \end{tabular}
\caption{Multigrid performance and comparison of the three precondtioners for MP1 with $k=160$, $256$ cells in the interior region, and $64$ cells in ECS-ABL on all four sides of the domain. ECS angle used is $\nicefrac{\pi}{6}$.}
\label{tab:exp1}
\end{table}

\begin{table}
 \begin{tabular}{|c|c|c|c|c|} \hline 
  \multirow{2}{*}{Preconditioner} & \multicolumn{2}{c|}{Multigrid} & mg cyc. & Bi-CGSTAB\\ \cline{2-3}
 & cyc, smooth., $\omega$ & mg-conv. , \# cycles & per prec. & iter , cputime \\ \hline
$CSL$ & F$_1^3$(1,1), $\omega$-Jacobi & 0.53 , 22 & \multirow{2}{*}{1} & \multirow{2}{*}{137 , 7m 34s}\\
$(\beta_1, \beta_2 = (-1, -0.4)$ & $\omega = 0.5$ & 13.4s & & \\ \hline
$CSG$ & F$_1^3$(1,1), $\omega$-Jacobi & 0.53 , 22 & \multirow{2}{*}{1} & \multirow{2}{*}{143 , 7m 36s}\\
$\theta_{\alpha} = \nicefrac{\pi}{17}$ & $\omega = 0.5$ & 14.4s & & \\ \hline
$QD$ & V(1,1), $\omega$-RB Jacobi & 0.15 , 8 & \multirow{2}{*}{1} & \multirow{2}{*}{357 , 19m 40s}\\
$\text{Re}(\lambda_0) = -16.88$ & $\omega = 1.0$ & 5.2s & & \\ \hline
 \end{tabular}
\caption{Multigrid performance and comparison of the three precondtioners for MP2 with $\lambda=7$, $k=4$. The domain is a square of 50 units. $512$ cells are used in the interior region, and $128$ cells in ECS-ABL on the east and the north side of the domain. ECS angle used is $\nicefrac{\pi}{6}$.}
\label{tab:exp2}
\end{table}

\begin{table}[!h]
 \begin{tabular}{|c|c|c|c|c|} \hline 
  \multirow{2}{*}{Preconditioner} & \multicolumn{2}{c|}{Multigrid} & mg cyc. & Bi-CGSTAB\\ \cline{2-3}
 & cyc, smooth., $\omega$ & mg-conv. , \# cycles & per prec. & iter , cputime \\ \hline
$CSL$ & F$_1^2$(1,1), $\omega$-Jacobi & 0.32 , 13 & \multirow{2}{*}{1} & \multirow{2}{*}{60 , 3m 9s}\\
$(\beta_1, \beta_2 = (-1, -0.6)$ & $\omega = 0.8$ & 6.45s & & \\ \hline
$CSG$ & F$_1^2$(1,1), $\omega$-Jacobi & 0.32 , 13 & \multirow{2}{*}{1} & \multirow{2}{*}{61 , 3m 10s}\\
$\theta_{\alpha} = \nicefrac{\pi}{13}$ & $\omega = 0.8$ & 6.3s & & \\ \hline
$QD$ & V(1,1), $\omega$-RB Jacobi & 0.17 , 8 & \multirow{2}{*}{1} & \multirow{2}{*}{164 , 9m}\\
$\text{Re}(\lambda_0) = -4.19$ & $\omega = 1.05$ & 1.2s & & \\ \hline
 \end{tabular}
\caption{Multigrid performance and comparison of the three preconditioners for MP3 with $k=2$. The domain is a square of 50 units. $512$ cells are used in the interior region, and $128$ cells in ECS-ABL on the east and the north sides of the domain. ECS angle used is $\nicefrac{\pi}{6}$.}
\label{tab:exp3}
\end{table}

\begin{table}
 \begin{tabular}{|c|c|c|c|c|} \hline 
  \multirow{2}{*}{Preconditioner} & \multicolumn{2}{c|}{Multigrid} & mg cyc. & Bi-CGSTAB\\ \cline{2-3}
 & cyc, smooth., $\omega$ & mg-conv. , \# cycles & per prec. & iter , cputime \\ \hline
$CSL$ & F$_1^4$(1,1), $\omega$-Jacobi & 0.32 , 13 & \multirow{2}{*}{1} & \multirow{2}{*}{210 , 18m 20s}\\
$(\beta_1, \beta_2 = (-1, -0.6)$ & $\omega = 0.8$ & 15.8s & & \\ \hline
$CSG$ & F$_1^3$(1,1), $\omega$-Jacobi & 0.31 , 12 & \multirow{2}{*}{1} & \multirow{2}{*}{160 , 14m 14s}\\
$\theta_{\alpha} = \nicefrac{\pi}{13}$ & $\omega = 0.8$ & 14.6s & & \\ \hline
$QD$ & V(1,1), $\omega$-RB Jacobi & 0.13 , 7 & \multirow{2}{*}{1} & \multirow{2}{*}{545 , 46m 40s}\\
$\text{Re}(\lambda_0) = -16.18$ & $\omega = 1.05$ & 9.4s & & \\ \hline
 \end{tabular}
\caption{Multigrid performance and comparison of the three precondtioners for MP3 with $k=4$, The domain is a square of 75 units. $768$ cells are used in the interior region, and $128$ cells in ECS-ABL on the north and the east sides of the domain. ECS angle used is $\nicefrac{\pi}{6}$.}
\label{tab:exp4}
\end{table}

\begin{figure}
\centering
\subfigure[Solution computed in Exp. 1]{\psfig{figure=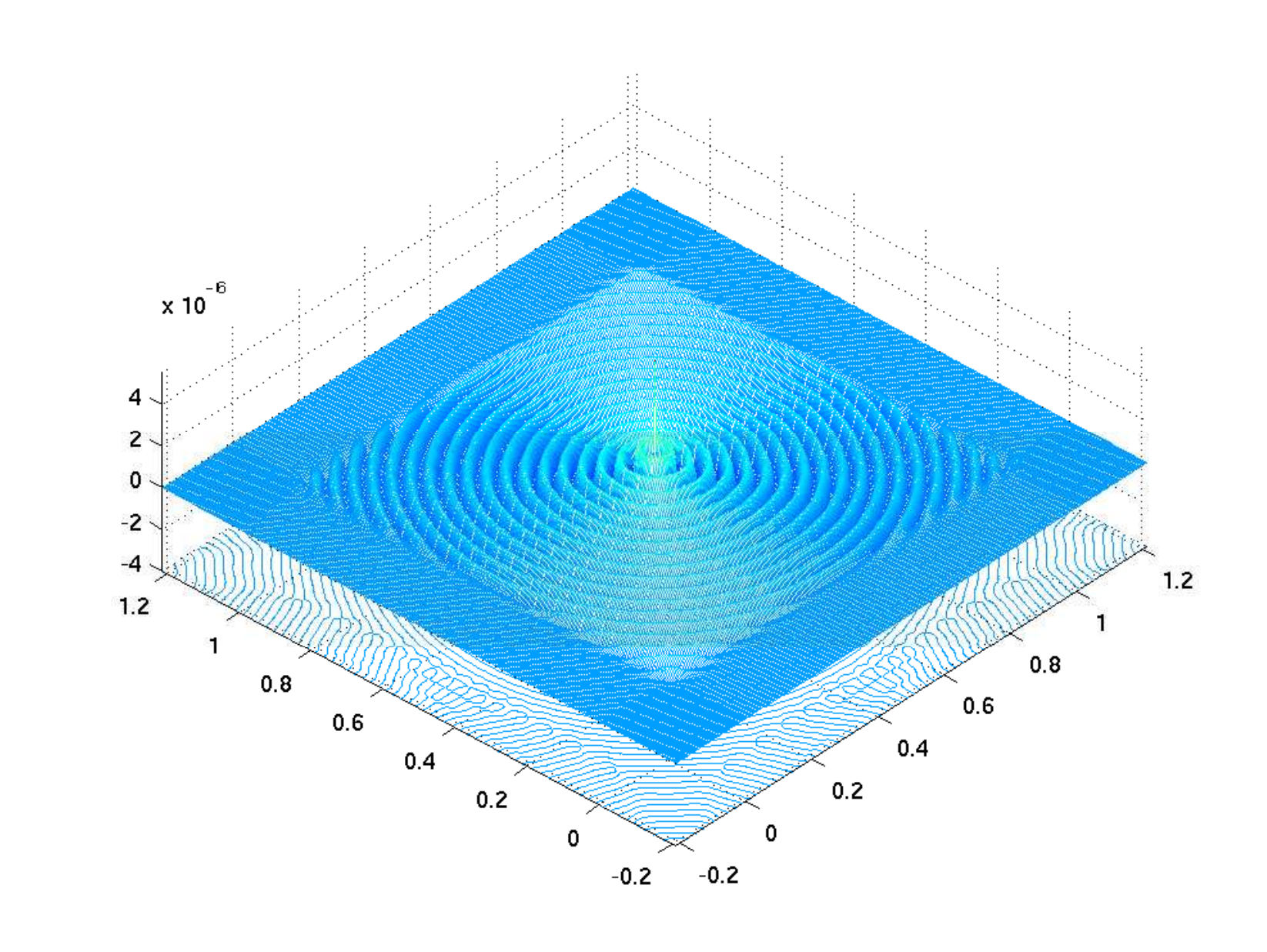, width=6cm}}
\subfigure[Solution computed in Exp. 2]{\psfig{figure=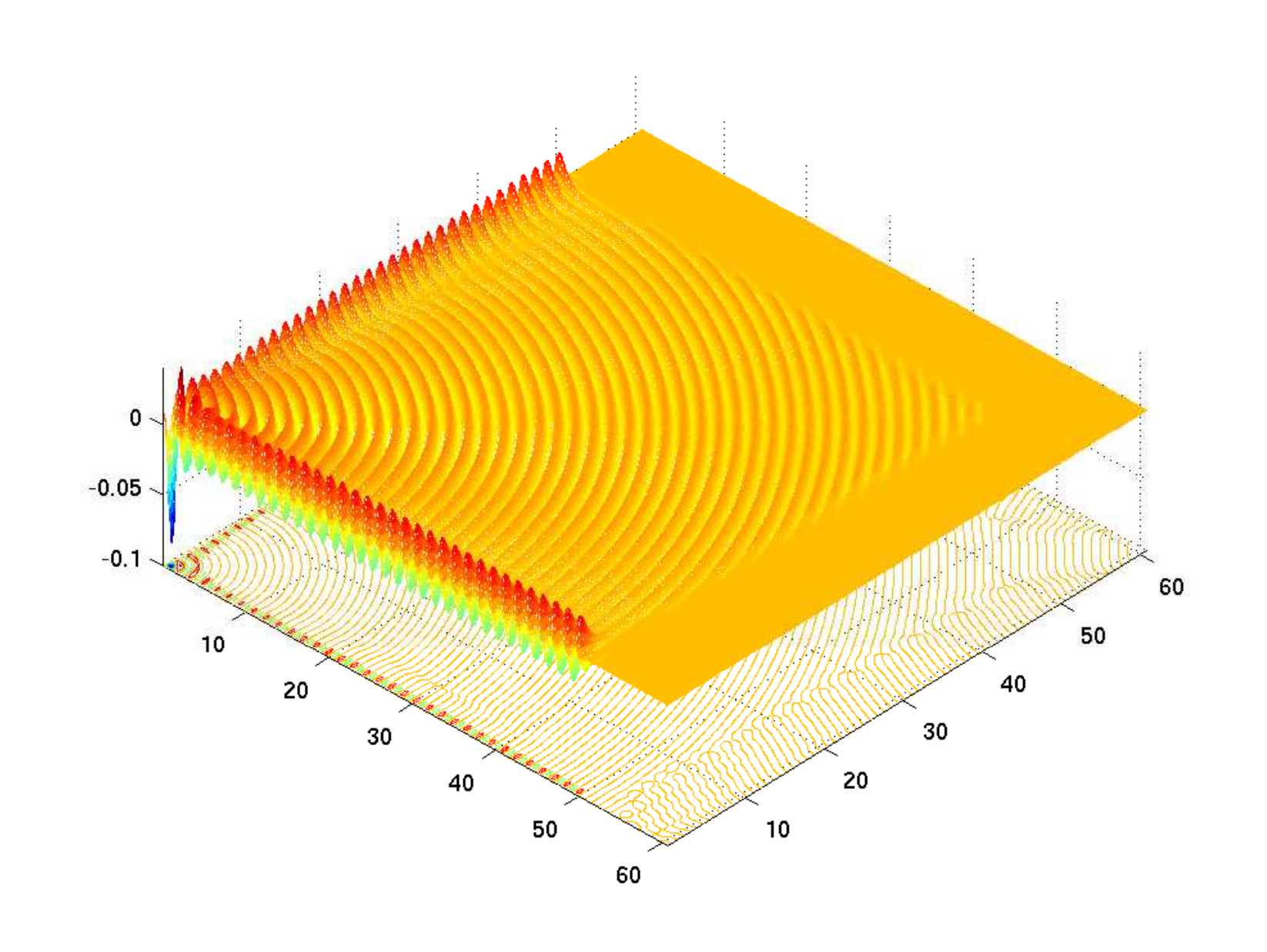, width=6cm}}\\ 
\subfigure[Solution computed in Exp. 3]{\psfig{figure=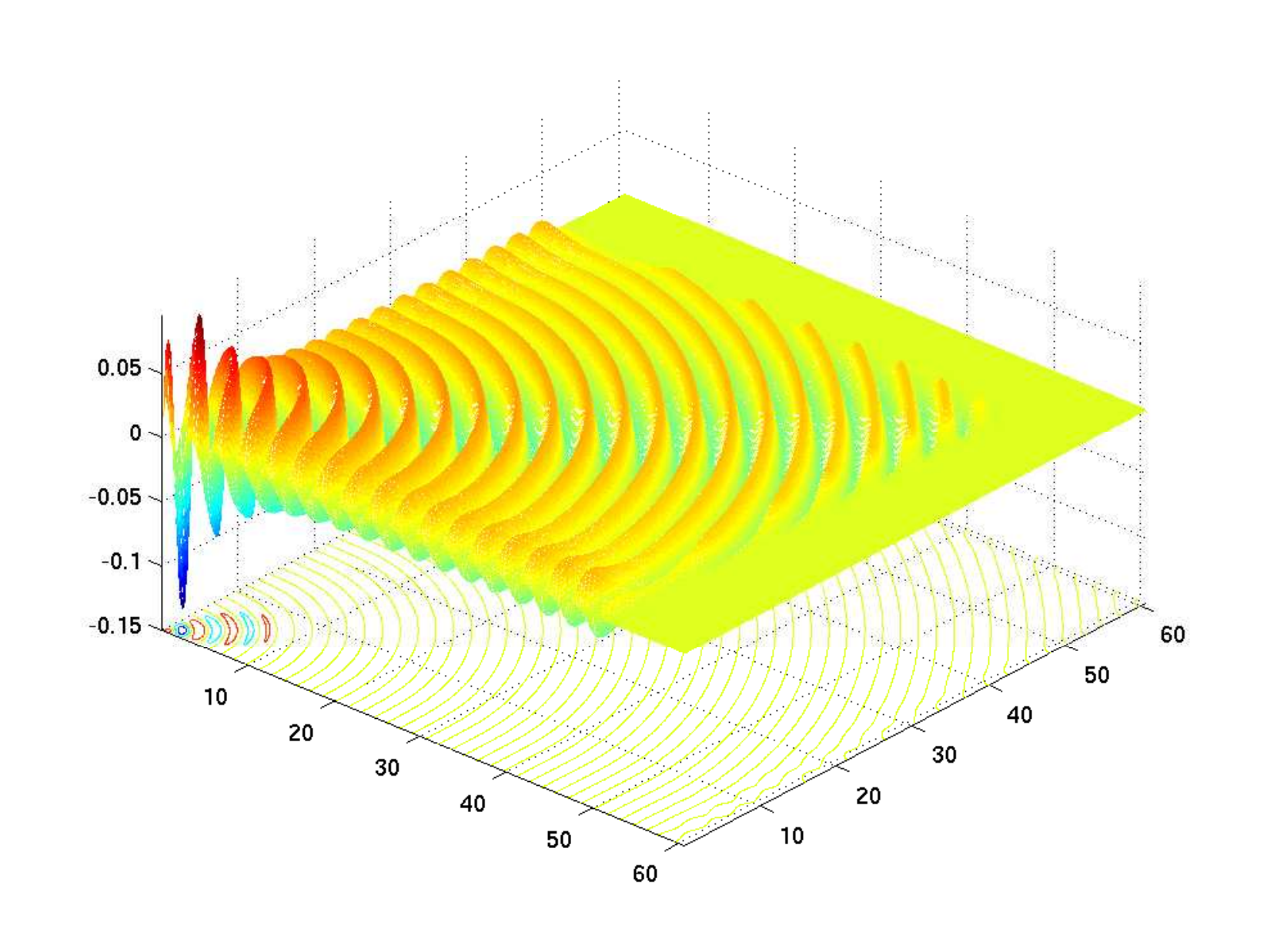, width=6cm}}
\subfigure[Solution computed in Exp. 4]{\psfig{figure=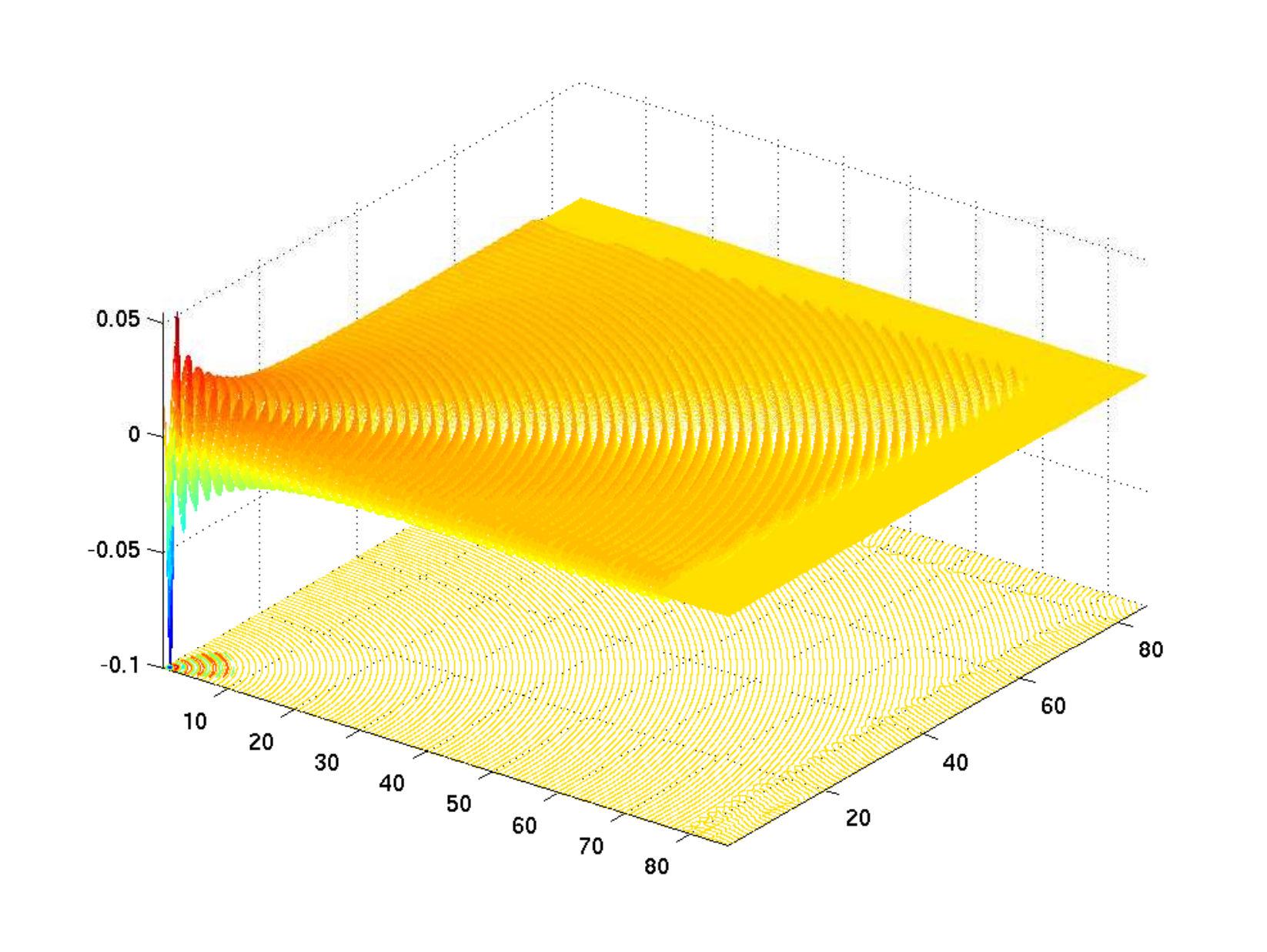, width=6cm}}
\caption{All these four solutions were computed for each of the 4 numerical experiments listed in the tables. The first solution as spherical waves ensuing out from the domain center. The other three solutions show evanescent waves, also known as single ionization, near the west and the south boundaries of the domain. At these edges the spatially dependent wavenumber grows exponentially in the model problems. (color online)}
\label{fig:solplots}
\end{figure}

\hbz{Numerical experiment results are reported for multigrid invertibility of the preconditioner and the observed efficiency of preconditioned Bi-CGSTAB. Multigrid invertibility is reported as the average multigrid convergence factor (mg-conv.) and the total number of cycles that the algorithm required to converge for the preconditioner taken as a stand-alone problem. mg-conv. is actually the geometric mean of the observed residual decay \edt{rates during multigrid cycles}, computed over the last $5$ cycles. The CPU-time is also reported. Bi-CGSTAB efficiency takes into account the number of iterations of the algorithm for convergence. Note that each Bi-CGSTAB iteration has two embedded multigrid cycles for preconditioning, i.e., one in each preconditioning step. The overall solution time is given as well.}

Results of the first experiment are listed in Table \ref{tab:exp1}. It is important to clarify that beating the $CSL$ or the $CSG$ precondtioner is not the aim of this work. We rather focus on obtaining a preconditioner which can come in close comparison to them in performance, and is comparatively much easier to invert. The $QD$ preconditioner takes around 3 times the number of iterations compared to the other choices. For this model problem (only) we also found that feeding in the preconditioner solution computed to a tolerance of \edt{10$^{-2}$}, to the Krylov method, as the starting guess, gives us a benefit of 50 iterations.

The rest of the experiments are listed in Tables \ref{tab:exp2},\ref{tab:exp3},\ref{tab:exp4}. They depict that the $QD$ preconditioner's performance comes within a factor of 3 of the other preconditioners even with strong spatial dependence in the wavenumber. The tuning effort is also much less, in fact, the relaxation parameters, the smoothing and the grid transfer methods can all stay constant. \hbz{To date, the authors are not aware of any scientific method that minimizes the $CSL$ shift or the $CSG$ rotation angle for different problems, without extra overhead. Note that these tunable parameters have a pivotal role in establishing $CSL$/$CSG$ superiority in speed over the $QD$ preconditioner. The experimental tables show the best cases for the $CSL$/$CSG$ preconditioners after they were hand-tuned for these parameters.} This points to the possibility that the $QD$ preconditioner might be used in an automatic solver setting. We tested the $QD$ preconditioner against the Laplacian preconditioner \hbz{(which can also be used in an automatic setting)} and found the $QD$ to be much superior in performance.

For multidimensional Helmholtz operators (including the 2D operator), the critical eigenvalue $\lambda_0$ used in the $QD$ preconditioner may be obtained from a one-dimensional counterpart, as is done in the current experiments. For Helmholtz problems with piecewise constant wavenumber, the maximum discrete wavenumber value may be used as a rough approximation of $\lambda_0$. However, this is also apt to bias the $QD$ preconditioner spectrum more to the right than is really required.

\section{Conclusions and Outlook}
\label{sec:conclude}
In this paper we showed that the Schr\"odinger equation for ionization problems can be decomposed into a coupled Helmholtz problem. This diagonal blocks of this coupled system consists of two-dimensional and three-dimensional Helmholtz problems. We propose Helmholtz model problems from these diagonal blocks. The blocks have homogeneous dirichlet boundaries at one side and exterior complex scaling absorbing layers  (ECS-ABL) at the other side. Finite difference discretization (for non-uniform grids) results in a pitchfork-shaped spectrum which is largely distributed in the fourth quadrant, but also has some parts crossing over in the third. Another property is that the spectrum is rather close to the real axis, and discrete problems are thus very challenging to solve iteratively. We solved them iteratively using the preconditioned Bi-CGSTAB method and also presented the quadrant definite ($QD$) preconditioner, which we derive from a time integration scheme for the Schr\"odinger equation. \edt{We tried using GMRES and restarted GMRES but found that for the current problems these methods failed to reach their superlinear convergence phase.} As a gross estimate we rate the efficiency of this preconditioner between the $CSL$/$CSG$  preconditioners and the Laplacian preconditioner, and it has the added advantage of having a multigrid favorable spectrum, i.e., its spectrum lies entirely in the fourth quadrant. This preconditioner can potentially be used in a automatic Helmholtz solver.  \wv{The advantage of the $QD$ preconditioner is that \edt{it} can be \edt{built} from standard multigrid components and it can be \edt{implemented} matrix-free which significantly reduces the memory use}.

\wv{Although we have used a low order discretization  of the differential operators and a low order absorbing boundary conditions, we believe that calculations with higher order methods will lead to \edt{similar} conclusions on the performance of the iterative method.}

Helmholtz problems, from an iterative perspective, can roughly be categorized into two classes which can be defined according to the available computational resources. One, where storing matrix operators is a possibility, and the other, where an iterative solution might have to be worked out using vectors alone. In the first situation, ILU(0) smoothing, and the Galerkin coarse grid operator used in a V(0,1) cycle render a very attractive multigrid method for preconditioner inversion. 

However, for the other class, the situation is comparatively much worse. First, for all preconditioners with a spectrum that leads to an efficient Krylov-subspace convergence, there is no appropriate smoother for multigrid.  Second, we have to do with re-discretizing the Helmholtz operator on the coarse grid. This seems to work with non-standard F-cycles (with multiple coarse grid recursions), which are expensive. In future, we intend to investigate, how smoothing may be enhanced for matrix-free Helmholtz solution contexts, as well as how to bring multigrid down to work in V cycles for preconditioner inversion. 
\section*{Acknowledgments}
This research was funded partially by \textit{Fonds voor Wetenschappelijk Onderzoek (FWO Belgium)} projects G.0174.08 and 1.5.145.10, by the \textit{Universiteit Antwerpen}, and by the Institute of Business Administration, Karachi, Pakistan. We wish to thank the sponsors sincerely for their support.

\bibliographystyle{elsarticle-num}
\bibliography{zubair}

\begin{thebibliography}{10}
\expandafter\ifx\csname url\endcsname\relax
  \def\url#1{\texttt{#1}}\fi
\expandafter\ifx\csname urlprefix\endcsname\relax\def\urlprefix{URL }\fi
\expandafter\ifx\csname href\endcsname\relax
  \def\href#1#2{#2} \def\path#1{#1}\fi

\bibitem{Wim05}
W.~Vanroose, F.~Martin, T.~N. Rescigno, C.~W. McCurdy, Complete photo-induced
  breakup of the {H}$_2$ molecule as a probe of molecular electron correlation,
  Science 310 (2005) ~1787--1789.

\bibitem{taylor2002computational}
K.~Taylor, {Computational challenges in atomic, molecular and optical physics},
  Philosophical Transactions: Mathematical, Physical and Engineering Sciences
  (2002) 1135--1147.

\bibitem{S79}
B.~Simon, The definition of molecular resonance curves by the method of
  exterior complex scaling, Physics Letters A 71 (1979) 211.

\bibitem{moiseyev1998}
N.~Moiseyev, {Quantum theory of resonances: calculating energies, widths and
  cross-sections by complex scaling}, Physics Reports 302~(5) (1998) 211.

\bibitem{B94}
J.-P. B\'erenger, A perfectly matched layer for the absorption of
  electromagnetic waves, Journal of Computational Physics. 114 (1994) 185--200.

\bibitem{CW94}
W.~C. Chew, W.~H. Weedon, A 3d perfectly matched medium from modified
  {M}axwell's equations with stretched coordinates, Microwave and Optical
  Technology Letters 7~(13) (1994) 599--604.

\bibitem{reps2009}
B.~Reps, W.~Vanroose, H.~bin Zubair, On the indefinite {H}elmholtz equation:
  exterior complex scaled absorbing boundary layers, iterative analysis, and
  preconditioning, Journal of Computational Physics 229 (2010) 8384–--8405.
\newblock \href {http://dx.doi.org/10.1016/j.jcp.2010.07.022}
  {\path{doi:10.1016/j.jcp.2010.07.022}}.

\bibitem{antoine2008review}
X.~Antoine, A.~Arnold, C.~Besse, M.~Ehrhardt, A.~Schaedle, A review of
  transparent and artificial boundary conditions techniques for linear and
  nonlinear schr\"odinger equations, Communications in Computational Physics
  4~(4) (2008) 729--796.

\bibitem{Bayliss83}
A.~Bayliss, C.~I. Goldstein, E.~Turkel, An iterative method for the {H}elmholtz
  equation, Journal of Computational Physics 49 (1983) 443--457.

\bibitem{Bayliss85}
A.~Bayliss, C.~I. Goldstein, E.~Turkel, On accuracy conditions for the
  numerical computation of waves, Journal of Computational Physics 59 (1985)
  ~396--404.

\bibitem{Yogi04}
Y.~Erlangga, C.~Vuik, C.~Oosterlee, On a class of preconditioners for solving
  the {H}elmholtz equation, Applied Numerical Mathematics 50 (2004) 409--425.

\bibitem{Yogi06}
Y.~Erlangga, C.~Oosterlee, C.~Vuik, A novel multigrid based preconditioner for
  heterogeneous {H}elmholtz problems, SIAM Journal on Scientific Computing 27
  (2006) ~1471--1492.

\bibitem{Achi77}
A.~Brandt, Multi-level adaptive solutions to boundary-value problems,
  Mathematics of Computation 31 (1977) ~333--390.

\bibitem{Stu82}
K.~St\"uben, U.~Trottenberg, Multigrid Methods: fundamental algorithms, model
  problem analysis and applications, Springer Berlin, 1982.

\bibitem{Trot01}
U.~Trottenberg, C.~Oosterlee, A.~Sch\"uller, Multigrid, Academic Press, 2001.

\bibitem{EEL01}
O.~C.~E. H.~C.~Elman, D.~P. O'Leary, A multigrid method enhanced by {K}rylov
  subspace iteration for discrete {H}elmholtz equations, SIAM Journal on
  Scientific Computing 23 (2001) 1291--1315.

\bibitem{meerbergen2009}
K.~Meerbergen, J.~Coyette, {Connection and comparison between frequency shift
  time integration and a spectral transformation preconditioner}, Numerical
  Linear Algebra with Applications 16~(1).

\bibitem{rescigno1999}
T.~Rescigno, M.~Baertschy, W.~Isaacs, C.~McCurdy, {Collisional breakup in a
  quantum system of three charged particles}, Science 286~(5449) (1999) 2474.

\bibitem{baertschy2001}
M.~Baertschy, T.~Rescigno, W.~Isaacs, X.~Li, C.~McCurdy, {Electron-impact
  ionization of atomic hydrogen}, Physical Review A 63~(2) (2001) 22712.

\bibitem{vanroose2006double}
W.~Vanroose, D.~Horner, F.~Mart{\'\i}n, T.~Rescigno, C.~McCurdy, {Double
  photoionization of aligned molecular hydrogen}, Physical Review A 74~(5)
  (2006) 52702.

\bibitem{arfken}
G.~B. Arfken, H.~J. Weber., Mathematical Methods or Physicists, Academic Press,
  1995.

\bibitem{baertschy2001solution}
M.~Baertschy, X.~Li, {Solution of a three-body problem in quantum mechanics
  using sparse linear algebra on parallel computers}, in: Proceedings of the
  2001 ACM/IEEE conference on Supercomputing (CDROM), ACM, 2001, p.~47.

\bibitem{mccurdy2004theoretical}
C.~McCurdy, D.~Horner, T.~Rescigno, F.~Martin, {Theoretical treatment of double
  photoionization of helium using a B-spline implementation of exterior complex
  scaling}, Physical Review A 69~(3) (2004) 32707.

\bibitem{horner2007two}
D.~Horner, F.~Morales, T.~Rescigno, F.~Mart{\'\i}n, C.~McCurdy, {Two-photon
  double ionization of helium above and below the threshold for sequential
  ionization}, Physical Review A 76~(3) (2007) 30701.

\bibitem{R95}
C.~M. Rappaport, Perfectly matched absorbing boundary conditions based on
  anisotropic lossy mapping of space, IEEE Microwave and Guided Wave Letters
  5~(3) (1995) 90--92.

\bibitem{TC98}
F.~L. Teixeira, W.~C. Chew, General closed-form {PML} constitutive tensors to
  match arbitrary bianisotropic and dispersive linear media, IEEE Microwave and
  Guided Wave Letters 8~(6) (1998) 223--225.

\bibitem{NB78}
C.~A. Nicolaides, D.~R. Beck, The variational calculation of energies and
  widths of resonances, Physics Letters A 65 (1978) 11.

\bibitem{mccurdyTR2004}
C.~W. McCurdy, M.~Baertschy, T.~N. Rescigno,
  \href{http://stacks.iop.org/0953-4075/37/R137}{Solving the three-body
  {C}oulomb breakup problem using exterior complex scaling}, Journal of Physics
  B: Atomic, Molecular and Optical Physics 37~(17) (2004) R137--R187.
\newline\urlprefix\url{http://stacks.iop.org/0953-4075/37/R137}

\end{thebibliography}

\end{document}